\documentstyle[11pt]{article}
\def\reals{\hbox{\sl I\kern-.18em R \kern-.3em}}
\def\quats{\hbox{\sl I\kern-.18em H \kern-.3em}}
\def\ints{\hbox{\sl Z\kern-.4em Z \kern-.3em}}
\def\nats{\hbox{\sl I\kern-.16em N \kern.05em}}
\def\rats{\hbox{\sl Q \kern-.83em\vrule height.59em depth0em \kern.87em}}
\def\complexes{\hbox{\sl\kern.50em I\kern-.50em C \kern.05em}}
\input boxedeps.tex
\SetepsfEPSFSpecial
\HideDisplacementBoxes

\begin{document}

\

\centerline{\bf BRAIDED CHORD DIAGRAMS}

\

\centerline{Joan S. Birman}
{\it \centerline {Barnard College of Columbia University}
\centerline {e-mail: jb@math.columbia.edu}
\centerline {Department of Mathematics, Columbia University, New York,
NY 10027}}

\

\

\centerline {Rolland Trapp}
{\it \centerline {California State University at San Bernardino}
\centerline {e-mail: trapp@math.csusb.edu}
\centerline {Department of Mathematics, C. S. U., San Bernardino, CA 92407}}

\

To appear in J. Knot Theory and its Ramifications {\bf 7}, No. 1 (Feb.1998)

\

\

\begin{abstract}

\noindent The notion of a braided chord diagram is introduced and studied. An
equivalence relation is given which identifies all braidings of a fixed chord
diagram. It is shown that finite-type invariants are stratified by braid index
for knots which can be represented as closed 3-braids. Partial results
are obtained about spanning sets for the algebra of chord diagrams of braid index
3.

\

\end{abstract}

\

\

\noindent {\it The first author dedicates her work in this paper to Professor
James Van Buskirk of the University of Oregon.  She thanks Jim for his
friendship and his many kindnesses over the years. }

\section{Introduction}
\label{Introduction}

The motivation for this work was the classification 
of links which are closed
3-braids by Birman and Menasco in \cite{[BM]} and the 
classification of finite
type invariants for knots which have braid index 2 in \cite{[T]}.  
Our goal was to
introduce the notion of braid index of chord 
diagrams, and then to study
finite-type invariants restricted to closed 3-braids.  
In particular, knowing
exactly which closed 3-braids represent non-invertible 
knots, and having a
precise classification theorem for all knots which are 
closed 3-braids, we hoped to
find examples which would prove that finite-type 
invariants detect
non-invertibility. 

We call the conjecture that such examples exist the 
NIC (non-invertibility
conjecture). The essential difficulty in finding 
examples to prove the NIC is
that present methods of calculation do not allow 
one to study weight
systems which have order larger than 9  \cite{[BN1]}. A
similar problem has arisen in other conjectures in 
this area. Since it is
known that every quantum invariant gives rise to an 
infinite family of finite
type invariants, it had been conjectured \cite{[BN1]} 
that all finite-type
invariants of a knot come from classical Lie Algebras, i.e. 
types A,B,C,D.  That conjecture  was based upon
the known data, i.e. weight systems of order up to 9. 
However, P. Vogel
proved that conjecture to be false \cite{[Vo]} when  he 
found a
counterexample (of order 32). We conjecture that a similar 
situation  exists with regard to the NIC. 

We remark that it
has been proved by G. Kuperberg \cite{[Ku]} that if finite-type invariants
contain enough information to completely classify prime {\it unoriented} knot
types,  then finite-type invariants are also strong enough to classify {\it
oriented} knots. Thus the question we are studying is
intimately tied to the larger question of whether finite-type invariants are a
complete set of algebraic invariants of knots. 

This paper develops a method and contains partial 
results  toward the goal of
finding examples which could settle the NIC.  
We begin in $\S2$ by proving some
general results about braided chord diagrams and 
braid index.  We show how two
braidings of the same diagram are related (Theorem 2.3) and 
give an algorithm for computing
the braid index of a diagram (Theorem 3.2).  The notion of
amalgamating fans  makes the computation of 
braid index significantly easier
(see Corollary 3.5) than it would be without 
such a notion.  Proposition 3.6
gives a lower bound for the braid index of a 
diagram which (as will be shown in
$\S$5) essentially allows us to characterize 
chord diagrams of braid index 3.

A fundamental
question is whether the value of a finite-type invariant on a
closed $b$-braid be computed from its values on diagrams of braid index 
at most $b$.
This can't be true if one uses Kontsevich's universal Vassiliev invariant in
the form in which it was presented in \cite{[BN2]} for the computation, because
the universal Vassiliev invariant evaluated on the trefoil involves diagrams with
braid index greater than 2. Thus we turned to the
actuality tables of \cite{[Va]} and \cite{[B-L]}. We say that finite-type
invariants are {\it stratified} by braid index if their values on
$b$-braids can be computed from a completed actuality table, using only diagrams
of braid index at most $b$.  The case for
$b=2$ was completed in \cite{[T]}, where finite-type invariants restricted to
closed 2-braids were characterized.  The difficulty for higher braid index is
showing that the crossing changes necessary to compute
$v(K)$ can be done without increasing the braid index.  To show this when
$K$ is a closed 3-braid, we begin by showing in Theorem 4.1 that 3-braid representatives of
$n$-diagrams are essentially unique.
We then use what we call standard singular
knot representatives to fill in the actuality tables.  Since
3-braid representatives are essentially unique, this
choice is well-defined.  Finally we use the algebraic structure of the
braid group to prove that finite-type invariants on closed 3-braids are
stratified by braid index (Theorem 4.4).  Thus, by the end of $\S$4 we have
shown that our program is sound: it is possible to develop a completed
actuality table for the computation of Vassiliev invariants of arbitrary order,
which computes such invariants for knots that are closed 3-braids and which
uses only singular knots and $n$-diagrams of braid index 3.  Since the chief
difficulty in collecting numerical data on Vassiliev invariants is that the
number of unknowns and the number of linear equations which one must solve grows
unmanageably large as one increases $n$ \cite{[BN1]}, and since both
the number of diagrams and the number of equations is necessarily
much smaller if one restricts to braid index
$\leq$3, it then follows that in principle one should be able to compute
higher order invariants in the 3-braid case than in the general case. 

Unfortunately, our results fell short of that goal.
We were able to determine all 4-term and 1-term relations in which all of the
chord diagrams have braid index 3, and also some more general linear relations
between them, but we do not at this time have a full understanding of all of the
linear relations which hold between chord diagrams of braid index 3. 
In $\S$5 we give strong forms of earlier theorems of Bar Natan \cite{[BN3]} and
also Lin \cite{[L]}, which take into account the braid forms of the 4-term
relations {\it and also} cyclic permutations of closed braided chord diagrams.
We give some examples which show why a better understanding of what we
call
$\lq\lq$stabilization" is necessary to complete the picture.

\section{Braidings of Diagrams}
\label{Braidings of Diagrams}

In this section we give a quick review of the background material. Then we
introduce the notion of braid index.  It is shown that every
diagram admits a braiding. We show that any two braidings of the same diagram
are related by a certain type of commutativity, cyclic permutation and
`stabilization'.

A {\it singular knot} is the image of a circle $S^1$ in $S^3$, under an 
immersion
$\iota:S^1\to S^3$ whose singularities are at most finitely many
transverse double points. Recall \cite{[B-L]} that any numeric knot
invariant
$v$ can be extended to singular knots via the Vassiliev skein relation:
$ v(K_p) =  v(K_{p_+}) - v(K_{p_-})$,
where $K_p$ is a singular knot with a transverse double point at $p$ and
$K_{p_+}$ (resp. $K_{p_-}$) is the singular knot with one less singular point
which is obtained from $K_p$ by `resolving' the singular point at $p$ to a
positive (resp. negative) crossing. All knot invariants will be considered to be
extended to singular knot invariants in this manner.  An invariant
$v$ is of {\it finite-type} if it vanishes on 
singular knots with more than $n$ double points, for some $n$.  It has
{\it order} $n$ if $n$ is the smallest such integer. It follows immediately that 
 an order-$n$ invariant evaluated on a singular knot $K^n$ with
$n$ double points depends only on the double points of the knot 
(because changing crossings in
$K^n$ changes the value of $v$ by its value on a singular knot with $n+1$ double
points, which is zero).  

A {\it chord diagram} is an $n$-diagram, for some $n$, i.e. an oriented circle
with $n$ chords. Examples are given in Figure 2.1. The numbers around the outer
circle label the chords, relative to an arbitrary starting point and an
anticlockwise orientation on the circle. We describe the diagram by a
{\it name}, i.e. the string of
$2n$ integers obtained by recording the labels of the chords as they are passed
while traversing the outer circle, as indicated in the examples in Figure 2.1.
Different choices of labeling and different starting points will yield
names which are cyclically equivalent after relabeling the chords. Given a
name, one may easily reconstruct the diagram $D$. Thus if two $n$-diagrams have
the same name, they are the same. An $n$-diagram can be thought of as the
addition of extra structure to the pair $(\iota,S^1)$, vis: think of $S^1$ as a
planar circle and of the endpoints of the chords as the points which
are identified under $\iota$.  A singular knot
$K^n$ which has a given
$n$-diagram
$D$ as its preimage is said to {\it respect} $D$.  

\ForceWidth{6in}
$$\EPSFbox{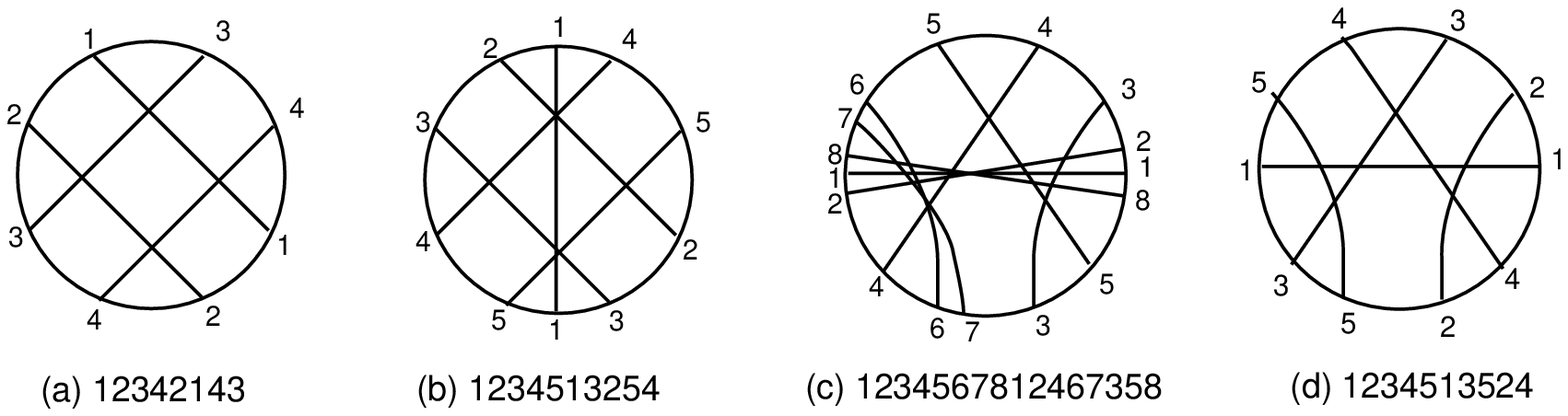}$$
\centerline{{\bf Figure 2.1: Examples of chord diagrams}}

\

If a knot invariant $v$ is of order $n$,
then $v(K^n)$ depends only on the diagram $D$ which $K^n$ respects.
The values which an invariant $v$ of order $n$ takes on $n$-diagrams are
determined by a system of linear equations whose unknowns are in one-to-one
correspondence with the collection of all possible 
$n$-diagrams. The linear relations in this system are known as the {\it 4-term} and 
{\it 1-term} relations \cite{[B-L]},\cite{[BN1]}.  A linear functional
$v$ on $n$-diagrams is called a {\it weight system} if it satisfies the 4 and 1-term relations. Kontsevich
\cite{[Ko]} proved that any weight system on $n$-diagrams extends to an invariant
of order $n$ on knots, and also that the extension is unique modulo invariants of
lesser order.  Thus the study of finite-type invariants can be reduced to the
study of weight systems on $n$-diagrams. Moreover, if the functional $v$ on
$n$-diagrams satisfies the four-term relation but not the 1-term relation, there
is a canonical way to modify
$v$ so that it satisfies the one-term relation as well \cite{[BN1]}.    

If we replace our planar circle $S^1$ by the 
closure of the standard m-braid representative
$\sigma_{m-1}^{-1}\sigma_{m-2}^{-1}\cdots\sigma_1^{-1}$ of the unknot, 
and choose the chords of
$D$ to be $n$ radial segments (relative to the braid axis) in this
representation, then we say that we have an {\it m-braid representative} or a
{\it braiding} of the
$n$-diagram. Cutting open the closed braid, we obtain a singular $m$-braid which
represents the diagram and has all of its chords as horizontal line segments. 
Two examples are given in Figures 2.2(a), where we have constructed two
braidings of the diagram $12342143$ of Figure 2.1(a). The associated open
braids are illustrated in Figure 2.2(b).  Clearly we may recover the closed 
braids by identifying the end points of the strands in the singular braid.

The chord algebra $C_m$ is defined as having generators $A(i,j)$
$1\le i < j \le m$ and relations

\[ [A(i,j),A(k,l)] = A(i,j)A(k,l)-A(k,l)A(i,j) = 0\]
\noindent
for $i,j,k,l$ distinct.  The geometric interpretation of a generator
$A(i,j)$ is a chord joining the $i^{th}$ and $j^{th}$ strand of the
identity $m$-braid, and multiplication in $C_m$ is concatenation.  While the
braids we are discussing are singular braids, we shall simply refer to them as
`braids'.

A word in the generators of 
$C_m$ can be closed in many ways to give a braided chord diagram, hence a
convention must be made.  If $W\in C_m$ is a monomial, then the {\it closure}
of $W$ is the braided chord diagram ${\bf W}$  obtained by adjoining the
standard
$m$-braid representative of the unknot to $W$ and forming the closure of the
resulting singular braid by identifying the top and bottom of each singular
braid strand. 
See Figure 2.2 for examples. The 4-diagram in Figure 2.2(a)
(which has the same name as
 that in Figure 2.1(a)) is the braided chord diagram associated to the closure
of the singular 4-braid $A(1,3)A(1,2)A(1,4)A(1,2)$ of Figure 2.2(b). The same
chord diagram is also defined by the closure of the singular 3-braid of Figure
2.2(b). 

\ 

$$\EPSFbox{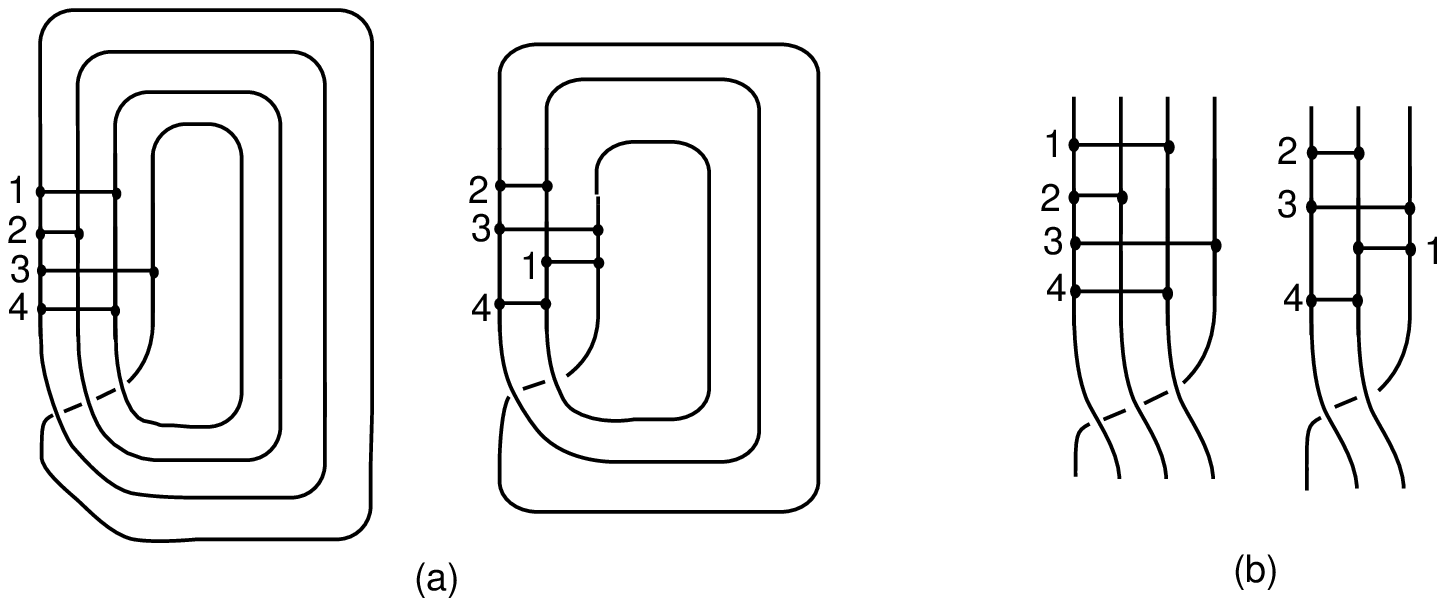}$$
\centerline{{\bf Figure 2.2: Braidings of chord diagrams}}

\

There is a natural injection
$i:C_l\rightarrow C_m$ for
$l\le m$ which acts as the identity on $C_l$.  Geometrically this injection is
realized by adding $m-l$ unused strands.  The reader should note that if $W\in
C_l$ then ${\bf W}$ and ${\bf i(W)}$ represent the same $n$-diagram since
the unused strands in $i(W)$ can be unkinked in ${\bf i(W)}$.  Thus we
adopt the convention that if $W$ is said to be in $C_m$, then $W$ contains
at least one generator of the form $A(i,m)$.  This convention is
unnecessary, but it simplifies some of the arguments which follow.

\vskip 6pt

\noindent
{\bf Lemma 2.1.} {\it  Any $n$-diagram has a braided chord diagram
representative.}

\vskip 6pt
\noindent
{\bf Proof.}  Let $N$ be a name for $D$, so that each integer
$r$ for $1\le r \le n$ occurs
twice in $N$.  Let $i_r$ be the position of the first occurrence of $r$
in $N$, and let $j_r$ be the position 
of the second.  We claim that the
word
$W= A(i_1,j_1)A(i_2,j_2)\cdots A(i_n,j_n)$ closes to a
representative of $D$.

Before proving this claim, we observe that $W\in C_{2n}$ and
every strand in $W$ has exactly one
endpoint on it.  These properties turn out to characterize braidings of
$D$ obtained in this way (see Lemma 2.2).

Now let's prove that ${\bf W}$ is a braiding of $D$.  To do so we simply
label the chord $A(i_r,j_r)$
of $W$ by $r$.  Now begin reading the name for ${\bf W}$ on the
left-most strand.  Note that by the way we labeled
the chords of $W$, the label on the $k^{th}$ strand is the $k^{th}$ entry
in $N$.  Thus $N$ is a name for ${\bf W}$ and ${\bf W}$ is a braid
representative for $D$. $\|$

\vskip 6pt

Note that the proof of the lemma gives an easy algorithm for constructing
a $2n$-braid representative from a name for an $n$-diagram.  This motivates
the following definition: Given a name $N$ for an $n$-diagram $D$, the braid
representative ${\bf W}_N$ of $D$ obtained as in the proof of Lemma 2.1
will be called the {\it canonical braiding} of $D$ relative to $N$.
Note that
in a canonical braiding the labeling of the chords determines the order of
the generators representing those chords.  Thus relabeling the chords
alters the canonical braiding by commutativity relations in $C_{2n}$.
Changing the starting point for reading $N$ alters the canonical
braiding by what will be called a {\it cyclic permutation}.

One further definition is needed.  Let $W$ be a word in the chord monoid
which closes to a braid representative of $D$.  We can obtain a name for
$D$ from $W$ by labeling the generators in $W$ in the order they occur, and
choosing the upper left-most strand as a starting point.
The resulting name for $D$ will be called the {\it standard name} for $D$
from $W$, and the process of obtaining the standard name from $W$ will
be called {\it reading a name} for $D$ from $W$.

Before showing how braidings of a
fixed diagram are related, we characterize  canonical braidings:

\vskip 6pt

\noindent
{\bf  Lemma 2.2} {\it  Let $D$ be an $n$-diagram and ${\bf W}$ a
braid representative of $D$.  Then ${\bf W}$ is a canonical braiding
of $D$ with respect to a name $N$ if and only if $W\in C_{2n}$ and
each strand of $W$ has exactly one endpoint on it.}

\vskip 6pt

\noindent
{\bf Proof.}  We have observed the only if portion of this
lemma in the proof of Lemma 2.1.  Conversely, let $W=
A(i_1,j_1)A(i_2,j_2)\cdots A(i_n,j_n)$ have the desired properties.
Let $N$ be the standard name for $D$ obtained by reading $W$.
Then the $k^{th}$ entry of $N$ is the label on the $k^{th}$
strand of ${\bf W}$, implying that ${\bf W}$ is the canonical
braiding of $D$ with respect to $N$. $\|$

\vskip 6pt

Recall that Markov's theorem implies that any two braid words with
isotopic closures are related by relations in $B_m$, conjugation in
$B_m$, and stabilization (where $B_m$ is the $m$-strand braid group
\cite{[B]}).  The corresponding theorem for braided chord diagrams is
similar.  The commutativity relations $[A(i,j),A(k,l)]=0$ in $C_m$ play
the role of the relations in $B_m$.  To interpret conjugation in the chord
diagram setting, consider the closure ${\bf W}$ of a word $W\in C_m$.
The final generator $A(i,j)$ in the word $W$ can be slid around the
closure portion of $\tilde{W}$ at the expense of increasing the indices
$i$ and $j$ by one.  Algebraically, we have that
the closed braid associated to  $W'A(i,j)$ is equal to the closed braid 
associated to $A(i+1,j+1)W'$ for
$1\le i < j < m$, also the same for $W'A(i,m)$ and $A(1,i+1)W'$.  Thus we
make one more definition: The operation of replacing the word $W'A(i,j)
\in C_m$ with $A(i+1,j+1)W'$ (with appropriate modification if
$j=m$) is called {\it cyclic permutation} in $C_m$.  The
operation which takes $C_m$ to
$C_{m+1}$  by replacing $WA(i,m)\in C_m$ with
$WA(i,m+1) \in C_{m+1}$, as in Figure 2.3, is called {\it stabilization}.  
We shall also use the term stabilization for the inverse of this operation.

$$\EPSFbox{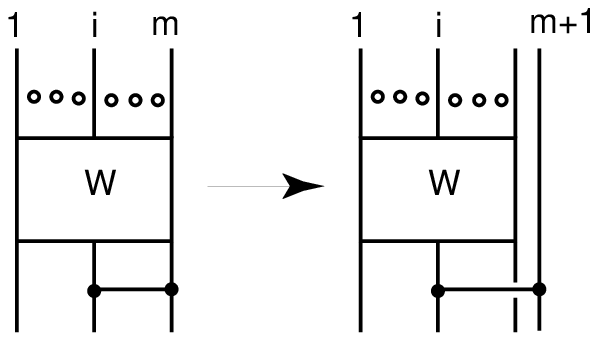}$$
\centerline{{\bf Figure 2.3: Stabilization}}

\

The notions of cyclic permutation and stabilization can be combined
to yield an operation that will be called generalized stabilization.
Generalized stabilization is adding or deleting a trivial loop anywhere
in the word $W$.  In the case of adding a trivial loop, this amounts to the
following algebraic operation.  Let $W= UA(i,j)V$
and replace it with $W'=U'A(i,j+1)V'$, where $U'$ is $U$ with all indices
greater than $j$ increased by one and $V'$ is $V$ with all indices $k\ge j$
increased by one.  This corresponds to
cyclically permuting $W$ until the generator $A(i,j)$ is last and has the
endpoint corresponding to $j$ on the final strand, stabilizing, and then
un-permuting to get $W'$. See Figure 2.4.  Deleting a trivial loop can
occur when $W=UA(i,j)V$ with $U$ void of the index $j$ and $V$ void of the
index $j-1$.  The word $W' = U'A(i,j-1)V'$ is again a braiding of the same
diagram as $W$, where $U'$, $V'$ are $U, V$ with indices greater than $j-1$
reduced by one.  Either of these processes is a generalized stabilization.
We will call adding a trivial loop an {\it increasing stabilization}, as
the braid index increases under this operation.  Similarly, deleting a trivial
loop will be called a {\it decreasing stabilization}.

$$\EPSFbox{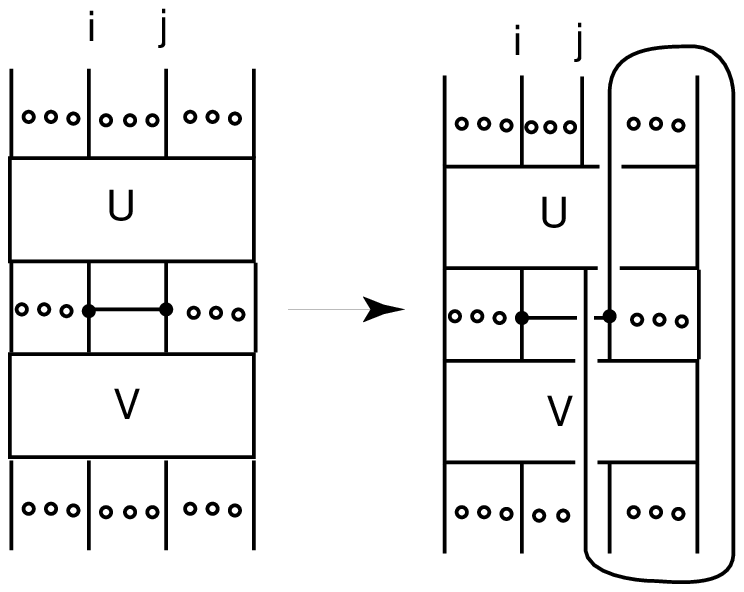}$$
\centerline{{\bf Figure 2.4: Generalized stabilization}}

\

Using these definitions, we describe an equivalence relation on braidings
which identifies all braidings of a fixed diagram:

\vskip 6pt

\noindent
{\bf Theorem 2.3.}  {\it Let ${\bf V}$ and ${\bf W}$ be two braidings of
a fixed diagram $D$.  The words $V$ and $W$ are related by a sequence of
the following moves:}
\begin{enumerate}
\item[(i)] {\it Commutativity relations in $C_m$,}
\item[(ii)] {\it  Cyclic permutation in $C_m$, and}
\item[(iii)] {\it  Stabilization.}
\end{enumerate}

\vskip 6pt

\noindent
{\bf Proof.}  We first show (Lemma 2.4) that any two canonical braidings are
equivalent under moves (i)-(iii) above. Then we show (Lemma 2.5) that these
moves can be used to change an arbitrary braid representative into a canonical
braiding. Lemmas 2.4 and 2.5 prove the theorem.

\vskip 6pt

\noindent
{\bf Lemma 2.4.} {\it  If ${\bf B}_N, {\bf B}_M$ are canonical braidings
of $D$ relative to the names $N$ and $M$, then ${\bf B}_N, {\bf B}_M$
are equivalent under moves of type (i) and (ii).}

\vskip 6pt

\noindent
{\bf Proof.}  In order to prove the lemma we recall how two names for
a diagram are related.  There are two choices to be made when
determining a name for a diagram; namely, the labels on the chords and
the choice of starting point for reading the labels.  Fixing the starting
point but changing the labels amounts to a reordering of the generators
in the canonical braid representative.  Since all the generators commute
in a canonical braid representative, this is accomplished by using
type (i) moves.  Fixing the labels and changing the starting point
changes the name by a cyclic permutation, which changes the canonical
braid representative by moves of type (ii).  This proves the lemma.
$\|$

\vskip 6pt

\noindent
{\bf Lemma 2.5} {\it  If ${\bf W}$ is a braiding of the diagram $D$ then
${\bf W}$ is equivalent under moves (i)-(iii) of Theorem 2.3 
to a canonical braiding for $D$.}

\vskip 6pt

\noindent
{\bf Proof.}  Let ${\bf W}$ be a braid representative for $D$, where
$W\in C_m$.  First we can assume that all the indices $1\le r\le m$ are
used in the word $W$, or that every strand of $W$ has an endpoint on it.
If not, then a finite number of decreasing stabilizations will reduce
$W$ to a braiding of $D$ without empty strands.  As decreasing stabilizations
are a consequence of moves (i) and (ii), we make this simplifying
assumption.

Let us assume, then, that every strand of $W$ has endpoints on it.  If $D$
is an $n$-diagram and each strand of $W$ has exactly one endpoint on it,
then $W\in C_{2n}$ and Lemma 2.2 implies that $W$ is a canonical
braid representative for $D$.  We are reduced to the case where $W$ has
at least one strand with more than one endpoint on it.  Note, however, that
if an increasing stabilization is performed just before the last endpoint
on a strand, then the number of endpoints on the given strand is reduced
by one and a new strand with a single endpoint is introduced. See Figure
2.5.  Thus a finite number of increasing stabilizations will create a
word $W'$ in which each strand has a single endpoint on it.  Thus $W'$
is a canonical braiding for $D$ by Lemma 2.2, proving the lemma and the
theorem.
$\|$

$$\EPSFbox{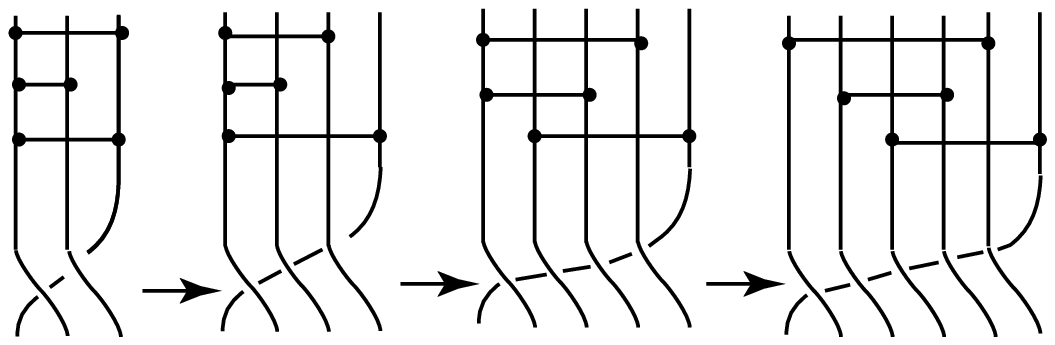}$$
\centerline{{\bf Figure 2.5}}

\

The reader who is familiar with \cite{[B-L]} or \cite{[BN1]} will notice that up
to this point we have not considered the 4 or 1-term relations at all. Theorem
2.3 does not consider equivalence modulo four-term relations, merely equivalence
of braidings of a single diagram.  The algebra
$$A^{hor} = \oplus_{m=0}^{\infty} C_m/\{4\ term\ relations\}$$
\noindent
is considerably more complex than $C_m$.  Since we have lots of work to do
before we get to $A^{hor}$ we defer our discussion of the 4 and 1-term 
relations to $\S$5. 
 
\section{Braid index}
\label{Braid index}

In this section we define and study the braid index of chord diagrams.
One reason braid index is of interest lies in the study of finite-type
invariants.  When studying weight systems on $n$-diagrams, restricting
to diagrams of a fixed braid index may simplify calculations.  We give
an algorithm for computing the braid index of a chord diagram, and
determine the effect of certain diagram characteristics on braid index.

\vskip 6pt

\noindent
{\bf Definition 3.1.}  {\it The braid index of an $n$-diagram $D$ is the least
number of strands used in any braid representative for $D$.}

\vskip 6pt

One might think that a more natural definition of braid index would
be the minimum braid index of any singular knot which respects $D$.
It turns out that this definition is equivalent to Definition 3.1.
Indeed, given a braid representation of $D$ of minimum braid index
one can obtain a singular knot which respects $D$ and has the same
braid index.  Merely  replace each chord $A(i,j)$ with the generator
$A_{ij}$ of the pure braid group with one of its crossings smashed
(see \cite{[B]} for a definition of the $A_{ij}$ and Figure 3.1).  The
singular knot which results from this operation will be called the {\it standard
knot} respecting ${\bf W}$.  This is well-defined as smashing
different crossings in the full twist of $A_{ij}$ results in singular
knots related by flypes.

$$\EPSFbox{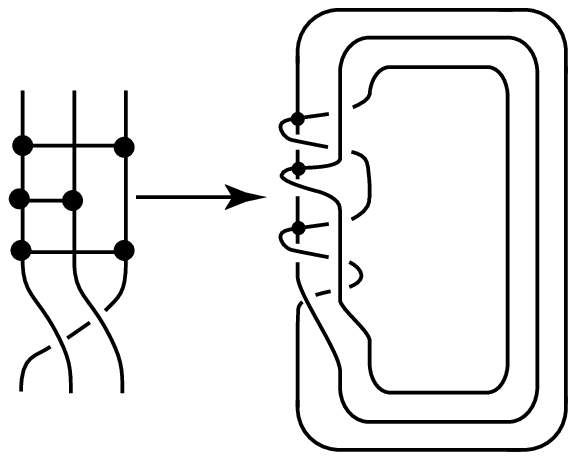}$$
\centerline{\bf Figure 3.1: The standard knot respecting the closure of a
singular braid}

\

Conversely,
given a braided singular knot $K$ which respects $D$ of
minimum braid index one
can construct a braid representative of $D$ with the same braid index.
To see this,
let $K$ have braid index $b$, and number the strands of $K$ by $1,\dots,b$.
To create a braid representative for $D$ from $K$ note that any braid
which closes to a knot must have an associated permutation that is a cycle.
Since every $b$-cycle is conjugate to (12$\dots b$), it must be possible
to cut the singular braid $K$ open somewhere and arrange it so that the
permutation on the cut-open braid is the $b$-cycle (12$\dots b$).  Let $\beta$
be the cut-open braid and $x$ the non-singular braid
$\sigma_{b-1}\sigma_{b-2}\dots\sigma_1$. Then $(\beta x^{-1})x$
closes to
$K$ and the associated permutation for
$\beta x^{-1}$ is the identity.  To finish constructing the braid
representative for $D$, begin with the identity
braid on $b$ strands.  Reading down the singular braid $K$, each time
a double point involves strands $i$ and $j$ connect strands $i$ and $j$
with a chord in the candidate for a braiding for $D$. See Figure 3.2.
After you pass completely through the singular braid, you have
the structure of the double points on the $b$-strand identity braid.
Closing the word in $C_b$ yields a braiding of $D$.
Thus the topological and algebraic definitions of braid index
agree.

\ForceWidth{6in}
$$\EPSFbox{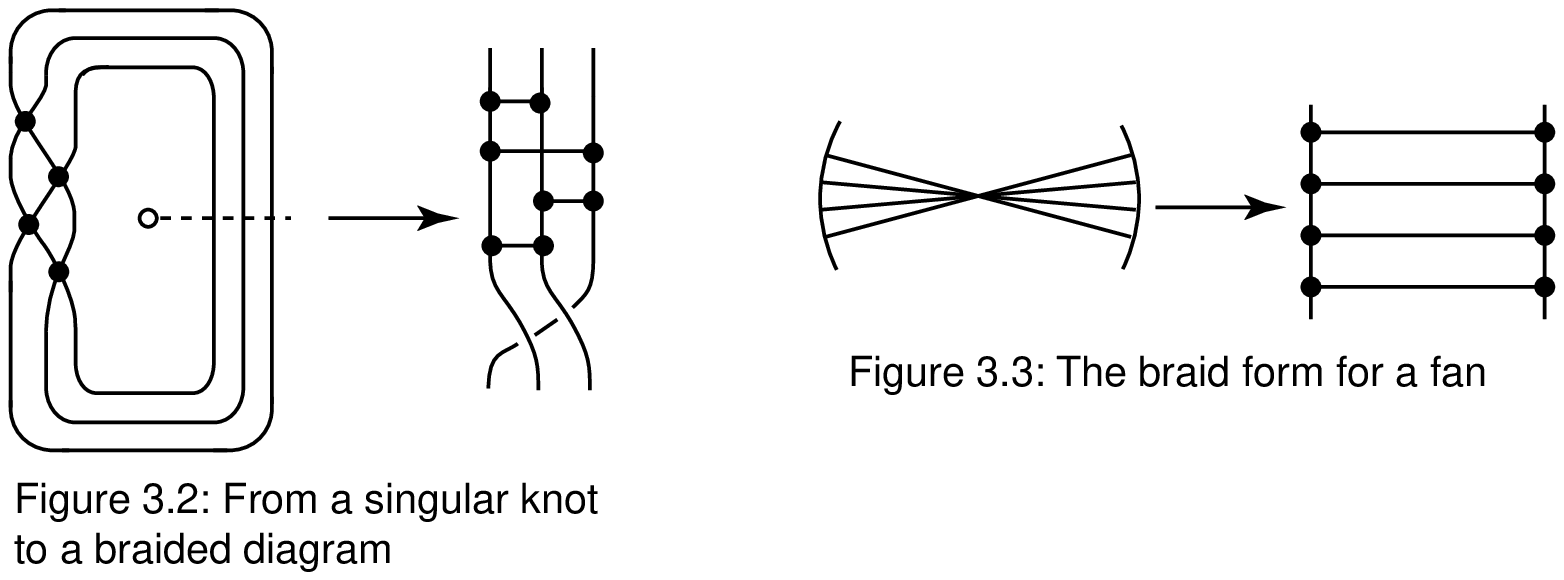}$$

We now give an algorithm for computing the braid index of an $n$-diagram
$D$.  To do so, we first recall the processes of increasing and
decreasing stabilization.  Recall that increasing stabilization is
the operation of replacing $W= UA(i,j)V$ with $W'=U'A(i,j+1)V'$,
where $U'$ is $U$ with all indices
greater than $j$ increased by one and $V'$ is $V$ with all indices $k\ge j$
increased by one.  The inverse of this process is a
decreasing stabilization.  Recall that one can easily recognize from a word
$W'=U'A(i,j)V'$ in the chord monoid $C_m$ when a braid index reducing
stabilization is possible.  If the index $j-1$ does not occur in $V'$ and
the index $j$ doesn't appear in $U'$, then a braid index
reducing stabilization is possible.  It has the effect of replacing the
word $W'=U'A(i,j)V'\in C_m$ with the word $W=UA(i,j-1)V \in C_{m-1}$, where
$U$ is $U'$ with all indices greater than $j$ reduced by one, and
$V$ is $V'$ with all indices at least as great as $j$ reduced by one
See Figure 2.4.  Decreasing stabilzations and
canonical braid representatives will be the essential ingredients of
our algorithm for computing the braid index.  The algorithm is as follows:

\vskip 6pt

\noindent
{\bf Step 1.}  Fix any chord of $D$, label it 1, and choose one endpoint
of chord 1, labeling it $\ast$.

\vskip 6pt

\noindent
{\bf Step 2.}  For each labeling of the remaining chords:
\begin{description}
\item[i.] Form the name resulting from beginning to read the labels at
$\ast$.
\item[ii.] Construct the canonical braiding of $D$ relative to that name.

\item[iii.]  Perform all possible decreasing stabilizations to the
canonical braiding, and record the braid index of the resulting braid.
\end{description}

\vskip 6pt
\noindent
{\bf Theorem 3.2.}  {\it The minimum of the braid indices found using the
above algorithm is the braid index of $D$.}

\vskip 6pt

\noindent
{\bf Proof.}  Let ${\bf W}$ be a braiding of $D$ of minimal braid index.
We will show that $W$ is cyclically equivalent to a word which is
obtained using only decreasing stabilizations from one of the names in
step 2 above.  As cyclic permutations don't change the braid index, this
implies that the minimum braid index is found in part (iii) of step 2.

Cyclically permute $W$ until the generator corresponding to chord 1 is
first and of the form $A(1,j_1)$, where the endpoint corresponding to $\ast$
is on the first strand.  Call the result $W$ again, and label the remaining
chords of $D$ in the order in which their corresponding generators
occur in $W$.  Now, as in the proof of Lemma 2.5, use
increasing stabilizations to achieve a canonical braiding for $D$ relative
to some name $N$.  Since increasing stabilizations do not change the
order in which generators occur, and since generators in canonical
braidings occur in the order given by the labeling of the chords,
we have the following.  The canonical braiding which $W$ increases to
is precisely the one found in part (ii) of step 2 using the ordering
of chords dictated by the order of generators in $W$. $\|$

\vskip 6pt

Once the first chord is chosen, there are $(n-1)!$ ways of labeling
the remaining chords.  Thus this algorithm grows factorially in
$n$.  Intuitively this algorithm is very simple.  Just place the
double points in a cycle, and connect corresponding double points
with edges always traveling counterclockwise.  The result is a closed
braid representative of $D$, so record its braid index.  There are
$(n-1)!$ cyclic orderings of $n$ points, and this algorithm just
checks each ordering. 

\vskip 6pt

We now discuss some techniques which help in determining braid index.
First we will discuss amalgamating fans, which allows one
to relate an $n$-diagram to a diagram with (possibly) fewer
chords and the same braid index.  Then we obtain a lower bound
on the braid index as being one more than the greatest number
of parallel chords of a diagram.  These techniques will be used
to characterize 2- and 3-braid diagrams. 

\vskip 6pt 

\noindent
{\bf Definition 3.3.} {\it Given a chord diagram $D$, a fan $F$ 
in $D$ is a collection of chords satisfying two
properties:}
\begin{itemize}
\item [(i)] {\it  Every chord in $F$ crosses every other chord in $F$,}

\item [(ii)]  {\it There are two arcs $\alpha_1$ and $\alpha_2$ on the outer
circle of $D$ such that each chord in $F$ has one endpoint in each
arc and no other chords in $D$ have endpoints in the $\alpha_i$.
See Figures 3.3 and 2.1 (c) for illustrations of fans.}
\end{itemize}

\

We remark that fans are easily recognized in names for $D$, they
result identical strings of integers occurring twice in the name.  A
maximal fan is a fan which is not contained in any other fan.  It is
clear that any chord in $D$ is contained in a unique maximal fan, possibly
consisting of the chord alone.  The process of replacing a fan $F$ of
chords in $D$ with a single chord, and labeling that chord with a weight
equal to the number of chords in $F$, will be called amalgamating the
fan $F$.  Let $D_a$ denote the diagram obtained from $D$
by amalgamating all maximal fans.  One retrieves $D$ from $D_a$ by simply
replacing each chord weighted $w$ with a fan of $w$ chords.  The reason
for discussing fans at this point is that fans do not affect the braid
index (see Proposition 3.4 and its corollary which follows).
Thus to compute the braid index of a diagram $D$ it suffices to compute
that of $D_a$, ignoring the weights on the chords.  When we speak of braid
representatives of $D_a$, then, we mean the unweighted diagram $D_a$.

\vskip 6pt

\noindent
\noindent
{\bf Proposition 3.4.}
\begin{itemize}  
\item [(i)] {\it  If $W_a = A(i_1,j_1)\dots A(i_n,j_n)$ 
closes to represent $D_a$, then the
word \\ $W=A(i_1,j_1)^{w_1}\dots A(i_n,j_n)^{w_n}$ closes to represent $D$.}
\item [(ii)]  {\it Any minimum braid index representative $W$ of $D$ is 
equivalent, using commutativity relations in $C_m$ and cyclic permutation, 
to the closure of a word \\ $A(i_1,j_1)^{w_1}\dots
A(i_n,j_n)^{w_n}$, where $W_a = A(i_1,j_1)\dots A(i_n,j_n)$ closes to
braid representative for $D_a$.}
\end{itemize}
\noindent
{\bf Proof.}  To prove (i) it suffices to check that raising a generator
$A(i_r,j_r)$ to a power $w_r$ has the effect of replacing the chord 
represented by $A(i_r,j_r)$ in $D$ by a fan of weight $w_r$.  

To see this we will compare the name $N_a$ obtained from 
$W_a$, that is $A(i_1,j_1)\dots A(i_n,j_n)$ to the name $N$ obtained from
$W = A(i_1,j_1)^{w_1}\dots A(i_n,j_n)^{w_n}$.  Increasing the power of a 
generator leads to parallel chords, and since the name $N$ is obtained
by reading down the strands in $W$, this leads to identical sequences
of integers (see the example in Figures 2.1 (c) and (d)).  Hence replacing
$A(i_r,j_r)$ in $W_a$ with
$A(i_r,j_r)^{w_r}$ has the effect of replacing both occurrences of $r$
in $N_a$ with a sequence of $w_r$ integers.  By the remark preceding
this proposition, this is tantamount to replacing chord $r$ in $D_a$ by
a fan of $w_r$ chords.

Now let us prove (ii).  Let ${\bf W}$ be a braiding of $D$, and choose a
starting point $\ast$ on the outer circle of $D$ which is not interior to
any fan.  Label the first chord after $\ast$ by 1, and cyclically permute
$W$ until the first generator corresponds to chord 1.  Further permute
until the first generator is $A(1,j_1)$, and the endpoint on the first 
strand corresponds to the endpoint of 1 which occurs immediately after
$\ast$ on $D$.  We will call this cyclic permutation $W$ as well. (In other
words, we permute $W$ until the name of $D$ obtained by reading $W$ 
is the same as if we start at $\ast$).  We would
like to show that using commutativity relations in $C_m$ alters $W$ to a 
word of the form $A(1,j_1)^{w_1}\dots
A(i_n,j_n)^{w_n}$, where each maximal fan of $D$ is represented in $W$ by
the power of a generator.  If this is the case, then removing the exponents
is equivalent to amalgamating fans and $W_a=A(i_1,j_1)\dots A(i_n,j_n)$ is
a braiding of $D_a$.

First we will show that
each chord in a maximal fan is represented in $W$ by the same generator,
then that relations in $C_m$ can be used to put $W$ in the desired form.
Let $F$ be a fan in $D$, and $r$, $s$ adjacent chords in $F$.  We claim
that $r$ and $s$ must be represented by the same generator in $W$.  Since
$r$ and $s$ are adjacent chords in $F$, the name for $D$ is of the form 
$N=  \dots rs\dots rs\dots$.  If $r$ and $s$ are represented by
different generators $A(i_r,j_r)$ and $A(i_s,j_s)$, then $i_r\le i_s$ and
$j_r\le j_s$ (this follows from the fact that $N$ is obtained by reading
$W$).  Since $N=  \dots rs\dots rs\dots$, no endpoints can occur between the
first $r$ and first $s$ in $W$.  If $i_r < i_s$, then a reducing 
stabilization would be possible.  This contradicts the fact that $W$ is
of minimal braid index, hence $i_r = i_s$.  Similarly one sees that $j_r =
j_s$.  Since adjacent chords in $F$ must be represented by the same generator
in $W$, all chords in $F$ must be represented by the same generator.

Now, if all the generators representing chords in $F$ are not in the same
syllable, then other generators must be between them.  However, since $F$
corresponds to identical sequences of integers in $N$, the indices of the 
intermittent generators must be distinct from those corresponding to chords
in $F$.  Thus the commutativity relations in $C_m$ can be used to make the
portion of $W$ corresponding to $F$ into a single syllable. $\|$

\vskip 6pt
\noindent
{\bf Corollary 3.5.}  {\it The braid indices of $D$ and $D_a$ are the same.}

\vskip 6pt
\noindent
{\bf Proof.}  Let $b(D)$ and $b(D_a)$ be the respective braid indices.
Part (i) of the previous proposition implies that $b(D)\le b(D_a)$, while
part (ii) implies the reverse inequality. $\|$

\vskip 3pt
Thus to find the braid index of an $n$-diagram, one can first amalgamate
all maximal fans and work with $D_a$.  This technique will be used in
proving that 3-braid representatives of chord diagrams are unique up
to cyclic permutation, for the most part.

We now establish some natural bounds on the braid index of a diagram.
\vskip 6pt
\noindent
{\bf Proposition 3.6.} {\it  If the set of chords of the diagram D has a subset
of $p$ parallel chords, then $b(D) \ge p+1$.}

\vskip 6pt
\noindent
{\bf Proof.}  We show that if $D_p$ is a diagram with $p$ non-intersecting
chords then $b(D_p) = p+1$.  The proposition follows by noting that if $D_p$
is the subdiagram of $D$ consisting of the $p$ parallel chords, then
$b(D) \ge b(D_p)$.

If $p=1$ then $D_p$ consists of a single chord and $b(D_p) = 2$.  Now suppose 
that $b(D_p) = p+1$ for all diagrams with $p$ non-intersecting chords, and
let $D_{p+1}$ be a diagram with $p+1$ parallel chords.  Choose an outermost 
chord of $D_{p+1}$, i.e. a chord which has all other chords on one side of 
it and no endpoints on the other, and call it chord $p+1$.  Let $D_p$ be the
diagram obtained from $D_{p+1}$ by removing the chord $p+1$, and let ${\bf W}$
be a minimal braid representative for $D_{p+1}$.  Since one side of chord
$p+1$ has no endpoints on it, one can cyclically permute ${\bf W}$ until the
generator representing chord $p+1$ is of the form $A(1,2)$.  Moreover, by
removing the generator $A(1,2)$ one introduces a decreasing stabilization,
hence $b(D_{p+1}) = b(D_p) + 1$.  By the induction hypothesis, $b(D_p) = p+1$,
and we are done. $\|$

\vskip 6 pt

Figure 2.1(b) gives an example which shows that the lower bound of
Proposition 3.6 is not sharp. By Proposition 3.6 the braid index of the diagram
in Figure 2,1(b) is at least 3. It is easy to see that there is a 4-braid 
representative, so in fact it is
at most 4.  For  3-braids we shall see in
$\S$5 that, modulo the 4-term relations there is a simple way to
recognize diagrams of braid index 3. Indeed, by Proposition 5.3, 
the braid index
cannot be 3, so it is in fact 4. The bound obviously is sharp in the case of
2-braids. In spite of much effort we were unable to find a combinatorial
property of diagrams (such as that given by Proposition 3.6) which would allow
us to compute the braid index by inspection, in the general case. As a weak
substitute we mention that there is a trivial upper bound:
\vskip 6pt
\noindent
{\bf Lemma 3.7.}  {\it If $D$ is an $n$-diagram, then $b(D) \le n+1$.}

\vskip 6pt
\noindent
{\bf Proof.}  Think of $D$ as being built by adding $n$ chords to a circle.
Each time you add a chord, you increase the braid index by at most one.
Since the braid index of the circle is 1, $b(D) \le n+1$.  $\|$

\section{Closed 3-braids}
\label{Closed 3-braids}

In this section we address the question of whether or not finite-type
invariants are stratified by braid index.  We give an affirmative answer
in the case of closed 2 and 3 braids.  A more 
precise statement of the question we address follows:  Can one choose 
representatives for chord diagrams in such a way that if $K$ is a
(honest-to-goodness) knot of braid index $b$ and $v$ is any finite-type
invariant, then $v(K)$ is determined by its values on chord diagrams 
of braid index at most $b$?

If one restricts their attention solely to the universal Vassiliev invariant,
the answer is clearly negative.  Computations made in \cite{[BN2]} show that
the universal Vassiliev invariant of the trefoil includes diagrams of
braid index higher than 2.  Thus one has to change their point of view
and consider completed actuality tables as in Vassiliev's initial work [Va]
and \cite{[B-L]}.  Our goal is to choose a singular knot which respects each
chord diagram in such a way that finite-type invariants are stratified
by braid index in this sense.

Recall \cite{[B-L]} that if one is given a completed actuality
table for a Vassiliev invariant $v$ and wants to determine $v(K)$ for
some knot $K$, one uses the Vassiliev skein relation and crossing
changes.  More precisely, one changes crossings of $K$ until it is the
unknot and using skein relation writes $v(K)$ in terms
of the value of $v$ on the unknot and on singular knots with one
double point.   The value of $v$
on the unknot is known, and one changes crossings in the singular knots
until the value of $v$ on them is in terms of known values and values
on knots with more double points.  The process continues until you
are working with singular knots with $n$ double points, where $n$ is
the order of $v$.  At each level, the known values of $v$ are the values
of $v$ on the choices of singular knots respecting chord diagrams.
One can think of this process as writing the original knot $K$ as
a sum of singular knots with increasing double points, where the singular
knots in the sum are those on which the value of $v$ is known (see \cite{[S]}
 for
a development of Vassiliev's theory along these lines).
The question then is can we choose singular knots respecting chord diagrams
at each level so that for any closed $b$-braid $K$, this sum contains only
singular knots of index at most $b$.

To show this when
$K$ is a closed 3-braid, we first show that 3-braid representatives of
$n$-diagrams are essentially unique.  We then use the standard singular
knot representative of ${\bf W}$ to fill in the actuality tables (if
3-braid representatives are not unique, then one has to show that this
choice is well-defined).  Finally we use the algebraic structure of the
braid group to prove that finite-type invariants on closed 3-braids are
stratified by braid index.

Before beginning we make some notational simplifications.  The generators
for $C_3$ are $A(1,2)$, $A(2,3)$, and $A(1,3)$ which we denote by
$a$, $b$, and $c$ respectively. See Figure 4.1.

\ForceWidth{6in}
$$\EPSFbox{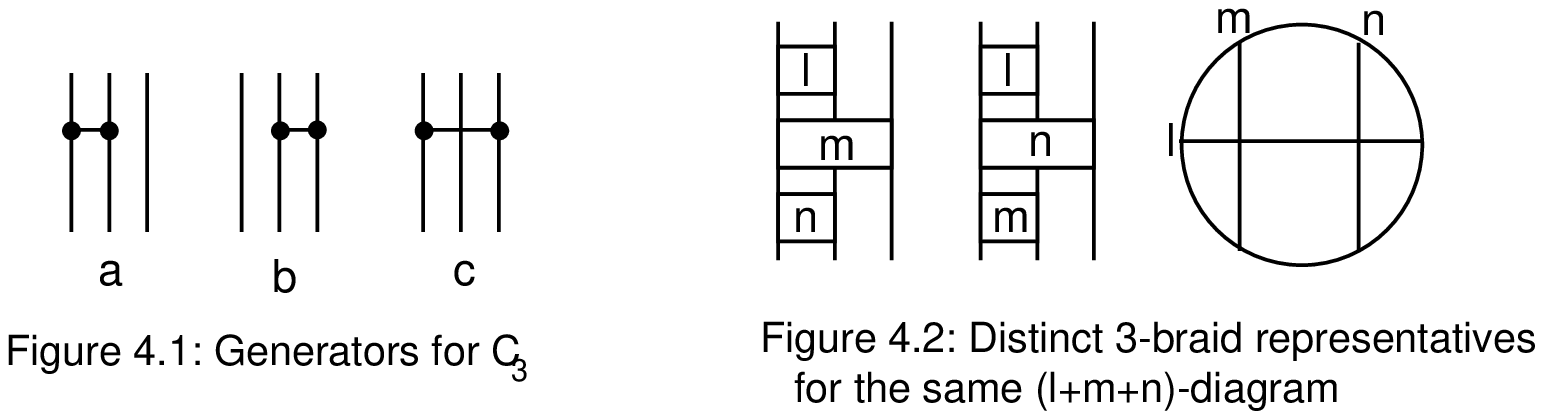}$$

Thus if $D$ is a 3-braid diagram, it can be represented by ${\bf W}$
where $W$ is a word in $a, b$, and $c$.  Recall that in braid index
considerations one can consider the amalgamated diagram $D_a$ instead.
Since all fans in $D_a$ consist of a single chord, the powers of generators
must be one in any word $W$ which closes to a representative of $D_a$.
In what follows we will usually consider the diagram
$D_a$, returning to our original diagram $D$ as necessary.  Note that in
$C_3$ all generators are cyclically equivalent.  In particular, under
cyclic permutation we have $a\rightarrow b\rightarrow c\rightarrow a$.

We now turn our attention to the question of uniqueness of 3-braid
representatives of chord diagrams.  We will show that the words $a^l
c^ma^n$, $a^lc^na^m\in C_3$ close to represent the same chord diagram,
but are not cyclically equivalent. See Figure 4.2.
We will also show that this is the
only case in which 3-braid representatives of a diagram are not unique
up to cyclic permutation.

\vskip 6pt

\noindent
{\bf Theorem 4.1.} {\it  Let $D$ be a 3-braid diagram, and let $V,W\in C_3$
be words which close to represent $D$.  Then either:}
\begin{enumerate}
\item[(i)] {\it $V$ and $W$ are cyclically equivalent, or}
\item[(ii)] {\it $V$ and $W$ are cyclically equivalent
to the words $a^lc^ma^n$, $a^lc^na^m\in C_3$.}
\end{enumerate}
\noindent
{\bf Proof.}  Let $D_a$ be the amalgamation of $D$.  As $V, W\in C_3$
close to represent $D$, Proposition 3.4 (ii) implies they are cyclically
equivalent to words in which all maximal fans appear as powers of generators. 
Thus we can cyclically permute $V, W$ until one obtains a braid representative
for $D_a$ from each simply by changing the exponents of generators to one.
We call the representatives of $D_a$ obtained in this way $V_a$ and $W_a$
respectively.  Note that cyclic permutations of $V_a$ and $W_a$ lift to
permutations of $V$ and $W$, so when we permute $V_a$ and $W_a$ we will
consider the words $V$ and $W$ permuted as well.

We would like some more structure, so we pick a labeling
$1,\dots,n$ of the chords of $D_a$ and a starting point $\ast$, which we
choose just before one of the endpoints of chord 1.  Let $N$ be the
name for $D_a$ obtained from this labeling and choice of starting point.
Now further permute $V_a$ and $W_a$ until the first generator in each
corresponds to chord 1 of $D_a$ and has the endpoint immediately following
$\ast$ on the first strand.  Thus $V_a$ and $W_a$ begin with an $a$ or
a $c$ generator, and reading either $V_a$ or $W_a$ yields
the name $N$.   There are now three cases to consider, one of which
immediately reduces to one of the other two.  Our goal in each case is to show
that
$V_a = W_a$, except in case (ii) of the theorem.  If $V_a = W_a$, it
follows that $V=W$.  Recall that $V$ and $W$ are obtained from $V_a, W_a$
by replacing single chords with fans (increasing the powers of exponents).
Given that $V_a=W_a$, the only way that we could have inequality
in $V$ and $W$ is if the same generator in $V_a,W_a$ corresponded to
different maximal fans in $D$.  This would mean that the same generator
in $V_a = W_a$ corresponded to different chords in $D_a$.  This can't
happen since reading either word yields the name $N$.

Once we show that $V_a = W_a$, then, we will know that $V$ and $W$
are cyclically equivalent.  The situation where $V_a \ne W_a$ happens
rarely, and we show that it happens only in situation (ii) of the
theorem.  Recall that $V_a$ and $W_a$ must have either an $a$ or a $c$
generator, and we study the following cases.
\vskip 3pt
\noindent
{\bf case 1:}  Both $V_a$ and $W_a$ begin with a $c$ generator.  

If this is the case, change the chosen endpoint $\ast$ of chord 1 and 
cyclically permute the words once.  You are then in case 2.
\vskip 3pt
\noindent
{\bf case 2:}  Both $V_a$ and $W_a$ begin with an $a$ generator.

We show that $V_a$ and $W_a$ are identical except in a few cases relying
heavily on the fact that reading either of the words gives the name
$N$ for $D_a$.

Since the first generator in $V_a,W_a$ is an $a$, and since the power of
each generator is 1, the second generator in each word is either a $b$ or
$c$.  If the second generators are identical we move on to the third, so
assume they are different.  Specifically, assume $V_a = acV'$ and $W_a=
abW'$, and that the $b$ generator in $W_a$ corresponds to chord $i$ in
$D_a$ while the $c$ in $V_a$ corresponds to chord $j$.  Reading $N$ from
$W_a$ we have $N= 1\dots 1i\dots i\dots$, and from $V_a$ we have
$N= 1j\dots 1\dots j\dots$.  Combining these partial names we see that
$j$ must follow the first 1 and $i$ the second in $N$.  Further, the 
second occurrence of $j$ in $N$ happens after the first $i$; therefore,
we have $N=1j\dots 1i \dots i\dots j\dots$ or $N=1j\dots 1i \dots j\dots 
i\dots$.  Lets consider each possibility separately.

Suppose that $N=1j\dots 1i \dots i\dots j\dots$, and recall that $V_a =
acV'$.  From $N$ we see that chord $i$ is parallel to both chords 1 and
$j$.  The fact that $i$ is parallel to 1 and its generator follows that
for 1 in $V_a$ implies that the generator corresponding to $i$ in $V_a$
must be a $b$.  However, a chord represented by a $b$ generator in $V'$
must intersect $j$, contradicting the fact that $i$ is parallel to both.
Hence this case cannot occur.

Now suppose $N=1j\dots 1i \dots j\dots i\dots$.  Since each generator in 
$V$ occurs to the first power, the third generator must be either an $a$ 
or a $b$.  If $V_a=acaV''$, and the second $a$ corresponded to the chord
$k$, then reading $V_a$ would give $N=1jk\dots 1k\dots j\dots$.  Hence this
case cannot occur and $V_a = acbV''$.  Moreover, name considerations imply
that the $b$ generator must correspond to the chord $i$.  This yields
that $N= 1j\dots 1i\dots ji\dots$.  Similar considerations with the word
$W_a$ yields $N= 1j\dots 1ij\dots i\dots$.  Taken together, these imply
that $N=1j\dots 1iji\dots$, and that $V_a = acbV''$ and $W_a = abaW''$
See Figure 4.3.

$$\EPSFbox{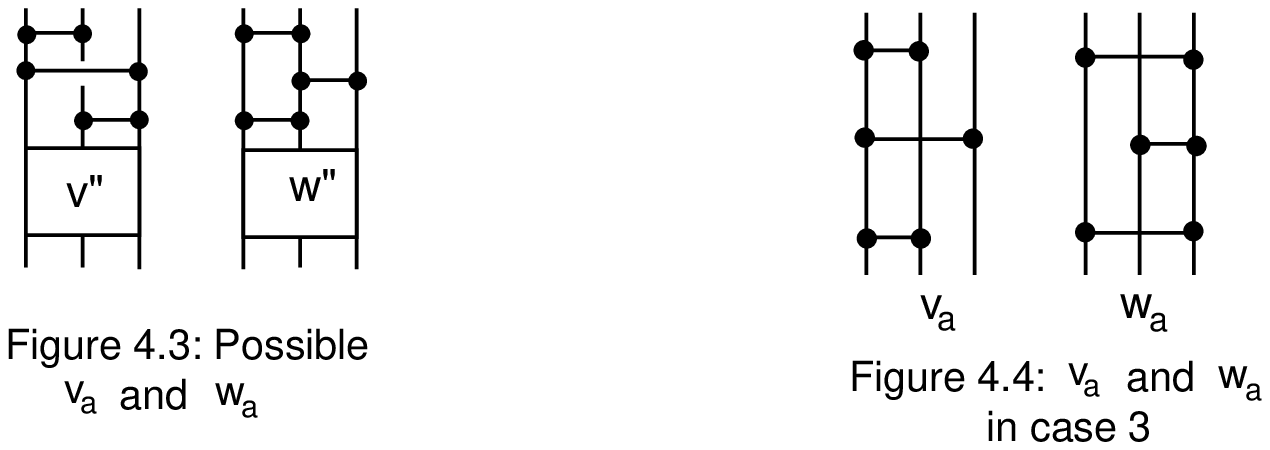}$$

\

Since no endpoints can occur between the first $i$ and second $j$,
the middle strand of $V''$ must be void of endpoints and $V''$
consists entirely of $c$ generators.  However, if $V''$ were non-empty
the $c$ generator would form a fan with the original $a$ generator in
$V_a$.  Since $V_a$ represents $D_a$ this cannot happen, and $V''$ is
empty.  Since $V_a$ and $W_a$ both represent $D_a$, it follows that
$W''$ is also empty.  

So far we have that if $V_a, W_a$ both begin with $a$ generators, but
have different second generators, then $V_a=acb$ and $W_a=aba$ (see Figure
4.2 for the diagram $D_a$).  The words $V_a=acb$ and $W_a=aba$ are 
cyclically equivalent; however, the cyclic permutation taking $W_a$
to $V_a$ does not preserve the initial point $\ast$ or the chords.  
(this arises because of the symmetry of the diagram $D_a$ in Figure
4.2).  However preservation of the point $\ast$ is essential if we consider 
the original diagram $D$ again.  If chords $1,i,j$ are weighted $n,m,l$
respectively, then we have $W= a^nb^ma^l$ and $V= a^nc^lb^m$.  These
are both representatives of the diagram in Figure 4.2, but are not
cyclically equivalent as the ordered triples of exponents are not.
These words, however, are cyclically equivalent to those in part (ii) of
the theorem, as it is easy to check.

Thus, if $V_a, W_a$ both begin with an $a$ generator and have different
second letters, then we are in case (ii) of the theorem.  We now have to
consider what happens if the first two generators of $V_a$ and $W_a$ are
identical, then the first three, etc.  Luckily we needn't argue as above,
but we prove a lemma and an easy induction takes care of this case.

\vskip 3pt
\noindent
{\bf Lemma 4.2.} {\it  Let $V_a$, $W_a$ be three-braid representatives of
$D_a$, cyclically permuted to start at the same point, and such that
$V_a = UV'$ and $W=UW'$ with $l(U) \ge 2$.  Then $V'$ and $W'$ are identical.}

\vskip 3pt
\noindent
{\bf Proof.}  Since $l(U)\ge 2$ at least two types of generators occur in 
$U$.  Suppose $U$ ends with an $a$ generator, then $N = u_1a\dots u_2a
\dots u_3\dots$ where the $u_i$ are the portions of $N$ coming from $U$
and the $\dots$ are the portions coming from $V'$ or $W'$.  We wish to 
show that the first generator in $V'$ and $W'$ are identical.  Since
$U$ ends in an $a$, $V'$ and $W'$ must begin with either a $b$ or a $c$.
If $V'$ begins with a $b$, then $N = u_1a\dots u_2 a k \dots u_3k\dots$ 
where $k$ is the chord represented by the first generator in $V'$.
If $W'$ begins with a $c$ generator, we have that $N= u_1al\dots u_2a
\dots u_3l\dots$, contradicting the fact that $k$ follows $u_3$.
Hence $W'$ begins with a $b$ as well.  A simple induction proves the
lemma in this case.  All other cases are similar. $\|$

\vskip 3pt
We can now consider 

\noindent
{\bf case 3:}  $V_a$ and $W_a$ begin with different generators.

Recall that neither $V_a$ nor $W_a$ can begin with a $b$ generator,
so assume that $V_a$ begins with an $a$ generator and $W_a$ with a $c$.
Hence $V_a = aV'$ and $W_a= cW'$.
Recall that the first generator in each word corresponds to chord 1 of
$D_a$, and reading each word gives the name $N$.  Thus $N= 1{\dots}_1 1
{\dots}_2$ where ${\dots}_1$ contains all endpoints on the first strand
of $V'$ and all those on the first and second strand of $W'$. See Figure
4.4.  This implies that no $b$ generators can occur in $V'$, otherwise 
we would have $N= 1\dots 1\dots i\dots i \dots$, and both $i$ endpoints
would have to be on the third strand of $W'$.  Since $W_a$ is a braiding
for $D_a$, this can't happen.  Similarly, one shows that there are no
$a$ generators in $W'$.  Thus we have $V_a = acaca\dots a$ and $W=
cbcbc\dots c$ ($V_a$ can't end in a $c$ generator, as that would
create a fan in $D_a$.  For the same reason $W_a$ ends in a $c$ generator).

Recall that the first generator in $V_a$ and $W_a$ corresponds to chord
1 in $D_a$.  The name obtained from reading $V_a$ is then 
$N= 1c_1a_1
\dots c_na_n1a_1a_2\dots a_nc_1c_2\dots c_n$ 
where the number $a_i$ 
(resp. $c_i$) is the label on the chord that the $i^{th}$ $a$ (resp. $c$)
generator in $V_a$ corresponds to.  Similarly, reading $W_a$ we have
$N= 1c_1c_2\dots c_nb_1\dots b_n1b_1c_1b_2c_2\dots b_nc_n$.  There is
no labeling of the $a_i, b_i, c_i$ which make these names identical
unless $n=1$ (note that the $c_i$ in $V_a$ could be different than the
$c_i$ in $W_a$).  For $n>1$, this contradicts the fact that reading
$V_a$ gives the same name for $D_a$ as reading $W_a$.  Thus the
only case in which $V_a$ and $W_a$ could start with different generators
is the case where $V_a = aca$ and $W_a = cbc$.

Again, even though these words are cyclically equivalent, the original
words $V$ and $W$ may not be.  Both $V_a$ and $W_a$ represent the diagram
$D_a$ of Figure 4.2, and we have already seen that if $V$ and $W$ are
not cyclically equivalent then we are in situation (ii) of the theorem.

Thus, in the case where both $V_a$ and $W_a$ begin with the same generator,
we have that 
\vskip 3pt
1.  $V_a = W_a$ and $V,W$ are cyclically equivalent, or

2.  $V_a, W_a$ are short and either cyclically equivalent or we are in case
(ii) of the theorem.
\vskip 3pt

In the case where $V_a$ and $W_a$ begin with different generators we've 
seen that $V_a, W_a$ are short and either cyclically equivalent or
examples of part (ii) of the theorem.

\vskip 3pt

The proof is complete.$\|$

\vskip 12pt

Now that 3-braid representatives of chord diagrams are essentially unique,
we have a well-defined method for choosing singular knots which respect
3-braid diagrams.  Recall that the standard knot respecting a braided
diagram ${\bf W}$ is the one obtained from the word $W$ by replacing
each $A(i,j)$ with the corresponding generator $A_{ij}$ of the pure 
braid group, and smashing one of the crossings.  To fill an actuality
table, one must choose singular knots which respect all chord diagrams.
For 3-braid diagrams, choose the standard knot corresponding to 
a braiding for each diagram.  We have the following

\vskip 3pt
\noindent
{\bf Corollary 4.3.}  {\it The choice of standard knots to respect 3-braid
diagrams is well-defined.}

\vskip 3pt
\noindent
{\bf Proof.}  One must show that standard knots respecting different 
braidings of the same diagram are isotopic.  We do this by picture.
First, cyclically equivalent braidings yield isotopic knots, as seen
in Figure 4.5.

$$\EPSFbox{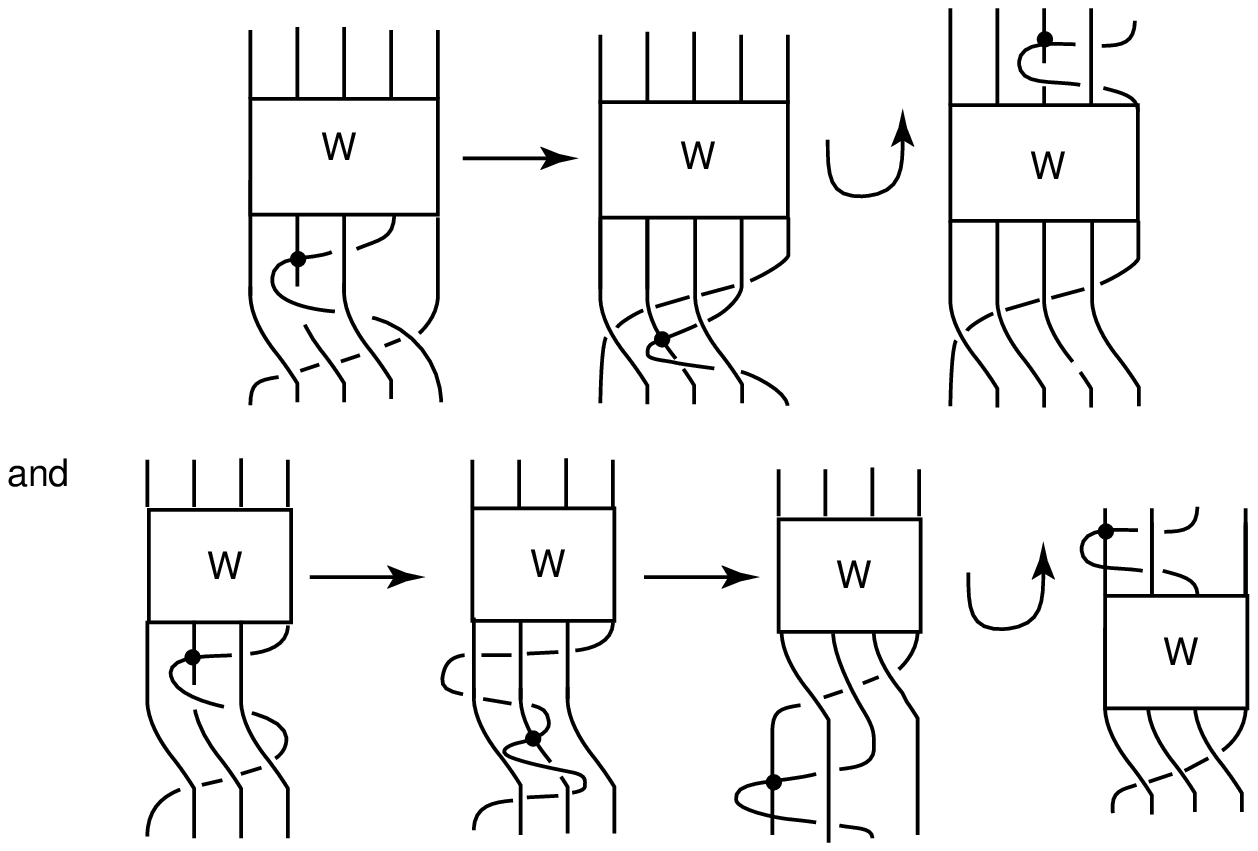}$$
\centerline{{\bf Figure 4.5: Isotopic standard knots}}

\

Secondly, one shows that the standard knots coming from the words
$a^lc^ma^n$ and $a^lc^na^m$ are isotopic via a braid-preserving flype
(note the similarity with the case of links in \cite{[BM]}.  See Figure 4.6.
$\|$

\ForceWidth{6in}
$$\EPSFbox{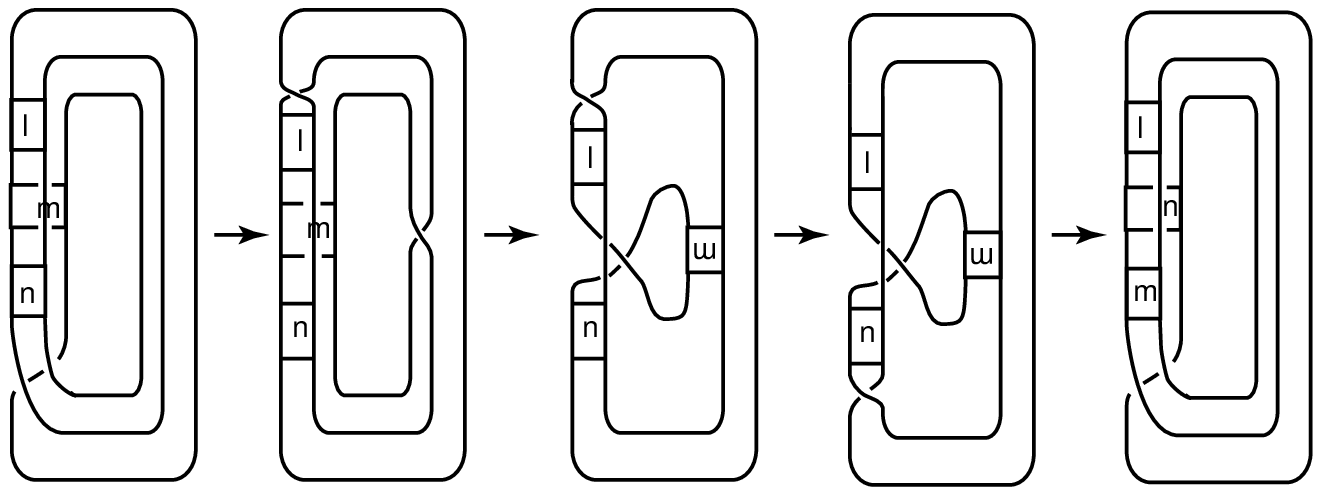}$$
\centerline{{\bf Figure 4.6: A braid-preserving flype}}

\

We can now prove the main theorem:

\vskip 6pt

\noindent
{\bf Theorem 4.4.} {\it  Finite-type invariants are stratified by braid index
for closed 3-braids.}

\vskip 6pt

\noindent
{\bf Proof.}  We must show that if we're given a closed 3-braid and
intend to evaluate a finite-type invariant on it, then we can perform
the necessary crossing changes without ever increasing the braid index.
It is well known \cite{[B]} that every element $\beta \in B_n$ can be
written uniquely in the form
$$\beta = \beta_2\beta_3\dots\beta_n\pi_{\beta},$$
where $\pi_{\beta}$ is a permutation braid and each $\beta_j$ belongs to
the subgroup of $B_n$ generated by the set $\{ A_{ij} : i<j\}$. We are 
concerned here
with the case $n=3$, and choose the permutation braid $\pi_{\beta} =
\sigma_2^{-1}\sigma_1^{-1}$.  Thus every knot which is a closed 3-braid can
be combed so that it's in the form $\beta_2\beta_3\sigma_2^{-1}\sigma_1^{-1}$.
Since the permutation braid is the one we've chosen to
close our braided chord diagrams, it is also the one appearing at
the end of every standard knot obtained from a braiding for a chord
diagram.  Thus switching crossings to undo the $\lq\lq$pure braid" part
will result in the standard knots we've chosen to respect our 3-braid
diagrams, and we are done. $\|$

\vskip 6pt

The difficulty with extending this proof to higher braid index is
that Corollary 4.3 no longer holds.  Then one can't be
assured that merely switching crossings that exist in the braid diagram
will yield the singular knot chosen to
respect a given chord diagram.  If one has to introduce other crossings,
there is no guarantee that braid index can be preserved.

\

\section{Four-term relations}
\label{Four-term relations}

In this section we make some concluding remarks about four-term relations in the
braid setting and discuss, briefly, the work which needs to be done to achieve
the goal we set forth at the beginning of this paper, i.e. to compute
finite-type invariants for knots which are closed 3-braids, using weight systems
which are restricted to chord diagrams of braid index 3.  

First, as is noted by many authors (e.g. \cite{[BN2]},
\cite{[L]}) the braid form of four-term relations is
$$ [A(i,j), A(j,k) + A(i,k)] = 0.$$
Using this relation it has been noted that you can comb elements in
$C_m/\{Four\ term\}$ to have the form shown in Figure 5.1 ( e.g. see
\cite{[BN3]}). That is, there is a spanning set for $C_m/\{Four\ term\}$ whose
elements have the special form $W_mW_{m-1}\dots W_2$, where $W_j$ is a word in
the generators
$A(1,j),A(2,j),\dots,A(j-1,j)$.

$$\EPSFbox{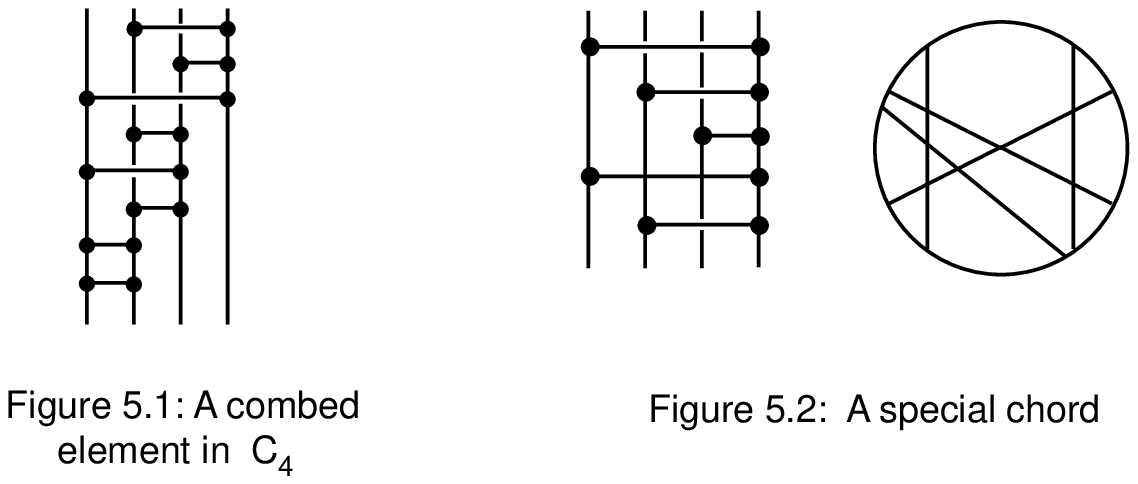}$$

We now show that there is a significant improvement in this normal
form when we allow cyclic permutations in addition to the four-term
relations. (i.e. ask what happens if we
$\lq\lq$close" words in $C_m$).  Letting $A_m = C_m/\{four\ term\ \&\ cyclic\
permutation\}$ we have the following result:

\vskip 6pt

\noindent
{\bf Proposition 5.1.} {\it  $A_m$ is spanned by words $W_m$ in the generators
$\{A(i,m), \ i=1,\dots,m-1 \}$ of $C_m$.}

\vskip 6pt

\noindent
{\bf Proof.}  Since we're working modulo four-term relations, we can
assume the words we begin with are already combed.  Now consider the final
generator in a combed word $W\in A$.  If $W = W'A(i,m)$ then we're done
since $W$ is combed.  Suppose that $W = W'A(i,j)$ with $j < m$, and let us
call each collection of generators in $W$ with the same last endpoint
a $\lq\lq$block" of $W$.  Then the last block in the combed $W$ must
involve the $j^{th}$ strand.

Cyclically permute $W$ to $A(i+1,j+1)W'$. This braid is not
combed.  We can comb $A(i+1,j+1)W'$ using, the four-term and
commutativity relations in $C_m$.  The result is a sum of combed words
$\Sigma W_i$ with the property that the final blocks of the $W_i$'s
have one less generator than the original $W$.  A finite number of
such cyclic permutations and combings (now with each $W_i$) replaces
$W$ with a sum of combed words in which the last block uses the $(j+1)^{st}$
strand.  If $j+1 = m$, we're done: otherwise, repeat the above procedure
until all words in the sum have the desired form. $\|$

\vskip 6pt

\noindent {\bf Remark 5.2:} Words in the generators $A(i,m)$ in
$C_m$ close to diagrams which have a nice property.  Before describing that
property, recall that in \cite{[C-D]} it is proved that, modulo the 4-term and
1-term relations, one can choose a basis of chord diagrams which have a chord
that intersects every other chord (a `special' chord).  The diagrams obtained by
closing
``one-block" words in $C_m$ either have a special chord or a single chord
could be added that intersects every other chord.  The key observation is that
if $W$ is a one-block word, then $A(1,m)W$ closes to a diagram with
a special chord.  This is clear since every chord in $W$ has exactly
one endpoint on the last strand, hence the endpoints are on either side
of the endpoints of the initial $A(1,m)$.  To see this, begin reading the
braid at the top of the $m^{th}$ strand. See Figure 5.2.  Thus
the above proposition is similar to the theorem of \cite{[C-D]},
yet not as precise. On the other hand, the one-term relation was not needed for
the proof of Proposition 5.1, where it is definitely needed to prove the
theorem of \cite{[C-D]}.

One might think that the reason that Proposition 5.1 is weaker than the
result of Chmutov and Dhuzhin \cite{[C-D]} is that we have
not used the 1-term relations, but more is at issue than that. Let's look at the
special case
$m=3$. Proposition 5.1 implies that $A_3$ is spanned by words in the generators
$b$ and $c$ for $C_3$  (see Figure 4.1).  If the word $W$ begins with
a $b$ generator, then the special chord has weight zero (i.e. the
word $cW$ represents a diagram with a special chord).  After amalgamating fans,
the cases where there is no special chord are thus the words $(bc)^k$.
Applying the methods used to prove the Main Theorem in \cite{[C-D]} to express
these special diagrams
$(bc)^k$ as a linear combination of diagrams with special chords (see Figure
5.3 for the cases $k=2,3,4$) we find that $k-1$ of the
terms on the right hand side involve diagrams with 3 parallel chords, which (by
Proposition 3.6) have braid index at least 4. Of course it is entirely
possible that there is some {\it other} linear combination of diagrams of braid
index 3 which will do the same job, we do not know. We note that the examples
in Figure 5.3 {\it do} use the 1-term relation. Clearly more work needs to be
done to understand this situation.

\ForceWidth{6in}
$$\EPSFbox{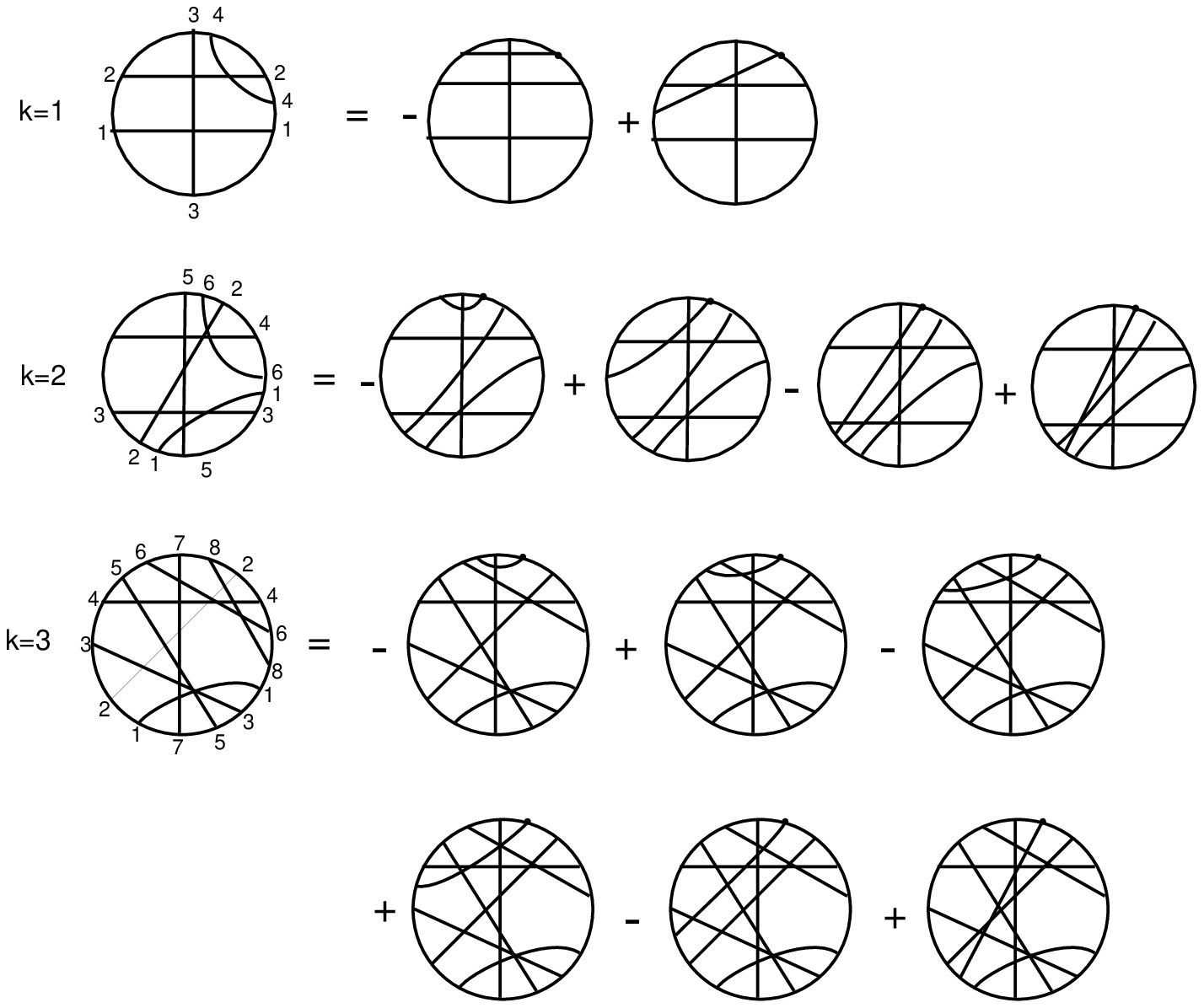}$$
\centerline{{\bf Figure 5.3}}

\

We continue our investigations. We next
show that 3-braid diagrams with special chords are easily recognized.
First note that a diagram with special chord (possibly of weight zero)
admits a particular type of name obtained as follows:  begin reading
at the special chord, numbering the chords as they are encountered.
The name so obtained is of the form $N = 123\dots n\sigma(1)\sigma(2)\dots
\sigma(n)$, where $n$ is the number of chords in $D$ and $\sigma$ is a
permutation on $n$ letters.  Given an diagram $D$ with special chord, call
the permutation $\sigma$ a permutation associated to $D$.  Note that in general
there may be several permutations associated to one diagram.  At the very
least, beginning to read at different ends of the special chord
yields the permutations $\sigma$ and $\sigma^{-1}$.  Recall that a descent
is a place where $\sigma(i) > \sigma(i+1)$.  We now prove

\vskip 6pt

\noindent
{\bf Proposition 5.3.} {\it  A diagram $D$ with special chord has braid index 3
if and only if an associated permutation for $D$ has a single descent.}

\vskip 6pt

\noindent
{\bf Proof.}  First suppose that $D$ is a 3-braid diagram, and let
${\bf W}$ be a 3-braid representative for $D$.  We must show that
there is a permutation $\sigma$ associated to $D$ with a single descent.
Cyclically permute $W$ until the generator corresponding to the special
chord is a $c$ generator at the beginning of the word, and call the
result $W$ again.  Thus $W = cW'$, and since $c$ intersects every other
chord $W'$ is a word in the generators $b$ and $c$.  Number the chords
by the order they appear in $W$ and begin reading the name for $D$ at the top
of the third strand.  Reading down the third strand yields the sequence
$123\dots n$, while continuing on the first strand yields an increasing
sequence consisting of the numbers of the $c$ generators.  Finally, reading
the second strand yields a second increasing sequence coming from the
$b$ generators.  Hence the name obtained in this manner has a single
descent.

Now suppose a name $N = 12\dots n\sigma(1)\sigma(2)\dots \sigma(n)$ for
$D$ has a single descent.  We construct a 3-braid representative for
$D$ as follows.  The sequence $\sigma(1)\dots\sigma(n)$ consists of two
increasing sequences.  Associate a $c$ generator to each number in the first
sequence and a $b$ generator to each in the second, then construct $W$ by
placing the $b$ and $c$ generators in numerical order.  Reading the word
$W$ beginning at the top of the third strand yields the original name $N$. $\|$

\vskip 6pt

\noindent {\bf Concluding Remarks:} From Remark 5.2 and the example given there
it is clear that a key obstacle to understanding the four-term relations in the
braid setting remains: to understand the interplay between the 4-term and
1-term relations and `stabilization'. We need to learn all of the relations
between chord diagrams of braid index 3 which are consequences of 4-term
relations in which chord diagrams of higher braid index appear. The spanning
set of Proposition 5.1, and that of \cite{[C-D]}, are too large to
be a basis for weight systems on chord diagrams of braid index 3, so unknown
additional relations clearly exist.

\

\noindent {\bf Acknowledgements}

\

\noindent The work of the first author was partially supported by NSF Grant DMS
94-02988; by the United States-Israel Binational Science
Foundation, under  Grant 95-00220/1; and by Barnard College, under a Senior
Faculty Research Leave, during the Spring Semester of 1997. She also
acknowledges the hospitality of the Mathematics Department of Technion
(Israel Institute of Technology).

\


\begin{thebibliography}{{\bf BM}9}
\markright{}
\baselineskip18pt

\bibitem{[B]} J.S. Birman, ``Braids, links and mapping class groups", {\it Ann.
of Math. Stud.}, {\bf 82} (1974), Princeton University Press, Princeton, N.J.

\bibitem{[B-L]} J.S. Birman and X.S. Lin, ``Knot polynomials and Vassiliev's
invariants", {\it Invent. Math.} {\bf 111} (1993), 225-270.

\bibitem{[BM]} J.S. Birman and W.W. Menasco, ``Studying links via closed braids
III: Classifying links which are closed 3-braids", {\it Pacific Journal of
Mathematics}, {\bf 161}, No 1 (1993), 25-113.

\bibitem{[BN1]} D. Bar-Natan, ``On the Vassiliev knot invariants", {\it Topology}
{\bf 34}, No. 2 (1995), pages 423-472.

\bibitem{[BN2]} D. Bar-Natan, ``Non-associative tangles", {\it Proc. 1993 Georgia
Int. Topology Conference}, to appear.

\bibitem{[BN3]} D. Bar-Natan, ``Vassiliev and quantum invariants of braids", {\it
Proc. Symp. in App. Math} {\bf 51} (1996), The Interface of Knots
and Physics, Editor L.H. Kauffman, 129-144.

\bibitem{[C-D]} S.V. Chmutov and S.V. Dhuzhin, ``An upper bound for the number of
Vassiliev knot invariants", {\it Jour. of Knot Theory and its Ramifications},
{\bf 3}, No. 2 (1994), 141-151.

\bibitem{[Ko]}  M. Kontsevich, ``Vassiliev's knot invariants", {\it Gelfand Seminar,
Part 2}, Advances in Soviet Math., {\bf 16} (1993), Amer. Math. Soc.,
Providence, Rhode Island, 137-150.

\bibitem{[Ku]} G. Kuperberg, ``Detecting knot invertibility", preprint,
U.C. Davis.

\bibitem{[L]} X.S. Lin, seminar, Columbia University, 1993.

\bibitem{[S]} T. Stanford, ``Finite-type invariants of knots, links, and
graphs", {\it Topology} {\bf 35} No. 4 (1996), 1027-1050.

\bibitem{[T]} R. Trapp, Twist sequences and Vassiliev invariants, {\it J.
Knot theory and its Ramifications}, {\bf 3} (1994), 391-405.

\bibitem{[Va]} V. A. Vassiliev, ``Cohomology of knot spaces", Theory of
singularities and its applications, (V.I. Arnold, ed.), Adv. Soviet Math.,
{\bf 1} (1990), Amer. Math. Soc., Providence, Rhode Island, 9-21.

\bibitem{[Vo]} P. Vogel, ``Algebraic structures on modules of diagrams",
preprint, August 1995. 


\end{thebibliography}
\end{document}